\input amstex 
\input amsppt.sty    
\magnification=\magstep1
\hsize 30pc 
\vsize 47pc 
\def\nmb#1#2{#2}         
\def\cit#1#2{\ifx#1!\cite{#2}\else#2\fi} 
\def\totoc{}             
\def\ign#1{}             
 
\redefine\o{\circ} 
\define\X{\frak X} 
\define\al{\alpha} 
\define\be{\beta}

\define\et{\eta}

\define\ka{\kappa}

\define\si{\sigma} 
 
\define\ph{\varphi} 
 
\define\ps{\psi} 
\define\om{\omega} 
\define\Ga{\Gamma}

\define\Om{\Omega} 
\redefine\i{^{-1}} 
\define\row#1#2#3{#1_{#2},\ldots,#1_{#3}} 
\define\x{\times} 
 
\define\Diff{\operatorname{Diff}} 
\define\Emb{\operatorname{Emb}} 
\define\dd#1{\frac{\partial}{\partial #1}} 
\define\g{{\frak g}} 
\define\ad{\operatorname{ad}} 
\def\today{\ifcase\month\or 
 January\or February\or March\or April\or May\or June\or 
 July\or August\or September\or October\or November\or December\fi 
 \space\number\day, \number\year} 
\topmatter 
\title  On the Geometry of the Virasoro-Bott group 
\endtitle 
\author  Peter W. Michor \\ 
Tudor S. Ratiu \endauthor 
\affil 
Erwin Schr\"odinger International Institute of Mathematical Physics,   
Wien, Austria  
\endaffil 
\address 
P\.W\. Michor: Institut f\"ur Mathematik, Universit\"at Wien, 
Strudlhofgasse 4, A-1090 Wien, Austria; and  
Erwin Schr\"odinger International Institute of Mathematical Physics,   
Pasteurgasse 6/7, A-1090 Wien, Austria. 
\endaddress
\email peter.michor\@esi.ac.at \endemail
\address
T\.S\. Ratiu: Department of Mathematics, University of California,
Santa Cruz, CA 95064, USA.   
\endaddress 
\email ratiu\@math.ucsc.edu \endemail 
\dedicatory \enddedicatory 
\date {\today} \enddate 
\thanks P.W.M\. was supported by `Fonds zur
F\"orderung der wissenschaftlichen                    
Forschung, Projekt P~10037~PHY'.
\endthanks
\thanks T.R\. acknowledges the partial support of NSF Grant DMS-9503273
 and DOE contract DE-FG03-95ER25245-A000.
\endthanks 
 \subjclass 58D05, 58F07, Secondary 35Q53 
\endsubjclass 
\abstract  We consider a natural 
Riemannian metric on the infinite dimensional manifold of all 
embeddings from a manifold into a Riemannian manifold, and derive its 
geodesic equation in the case $\Emb(\Bbb R,\Bbb R)$ which turns out 
to be Burgers' equation. Then we derive 
the geodesic equation, the curvature, and the Jacobi 
equation of a right invariant Riemannian metric on an  infinite 
dimensional Lie group, which we apply to $\Diff(\Bbb R)$, 
$\Diff(S^1)$, and the Virasoro-Bott group. Many of these results are 
well known, the emphasis is on conciseness and clarity.  
\endabstract 
\endtopmatter 
\leftheadtext{\smc Peter W. Michor, Tudor Ratiu} 
 
\document 
 
\heading Table of contents \endheading 
\noindent 1. Introduction \leaders \hbox to 
1em{\hss .\hss }\hfill {\eightrm 1}\par  
\noindent 2. The general setting and a motivating example 
\leaders \hbox to 1em{\hss .\hss }\hfill {\eightrm 2}\par  
\noindent 3. Right invariant Riemannian metrics on Lie groups 
\leaders \hbox to 1em{\hss .\hss }\hfill {\eightrm 4}\par  
\noindent 4. The diffeomorphism group of the circle revisited 
\leaders \hbox to 1em{\hss .\hss }\hfill {\eightrm 9}\par  
\noindent 5. The Virasoro-Bott group and the Korteweg-de 
Vries-equation \leaders \hbox to 
1em{\hss .\hss }\hfill {\eightrm 11}\par   

\head\totoc\nmb0{1}. Introduction \endhead

We consider a natural Riemannian metric on the infinite dimensional 
manifold of all embeddings from a manifold into a Riemannian 
manifold, derive its geodesic equation in the case 
$\Emb(\Bbb R,\Bbb R)$ which turns out to be Burgers' equation. For 
the general case see \cit!{9}.
Then we give a careful exposition of the derivation of the geodesic 
equation, the curvature, and the Jacobi equation of a right invariant 
Riemannian metric on an  infinite dimensional Lie group. This is a 
careful presentation and extension of results in \cit!{1}, \cit!{2}, 
\cit!{3}. The formulas obtained in this way are then applied to 
$\Diff(\Bbb R)$, $\Diff(S^1)$, and the Virasoro-Bott group, where the 
geodesic equation is the Korteweg-de~Vries equations. This is due 
to \cit!{8}, \cit!{10}, \cit!{23}, and also \cit!{22}. A fast 
overview on the geometry of the Virasoro-Bott group can also be found 
in \cit!{24}. The emphasis of this paper is on a unified setting for 
these results, and on conciseness and clarity. Thanks to Hermann 
Schichl and to the referee for detailed checks of the computations. 
 
\head\totoc\nmb0{2}. The general setting and a motivating example  
\endhead 
 
\subhead\nmb.{2.1}. The principal bundle of embeddings \endsubhead 
Let $M$ and $N$ be smooth finite dimensional manifolds, connected and  
second countable without boundary, such that $\dim M\leq\dim N$. 
The space $\Emb(M,N)$ of all embeddings (immersions which are 
homeomorphisms on their images) from $M$ into $N$ is an 
open submanifold of $C^\infty(M,N)$ which is stable under the right 
action of the diffeomorphism group of $M$. 
Here $C^\infty(M,N)$ is a smooth manifold modeled on spaces of 
sections with compact support $\Ga_c(f^*TN)$. In particular 
the tangent space at 
$f$ is canonically isomorphic to the space of vector fields along $f$ 
with compact support in $M$. If $f$ and $g$ differ on a non-compact 
set then they belong to different connected components of 
$C^\infty(M,N)$. See \cit!{19} and \cit!{14}.
Then $\Emb(M,N)$ is the 
total space of a smooth principal fiber bundle with structure 
group the diffeomorphism group of $M$; the base is called $B(M,N)$, it 
is a Hausdorff smooth manifold modeled on nuclear (LF)-spaces. 
It can be thought of as the "nonlinear Grassmannian" of all 
submanifolds of $N$ which are of type $M$. This result is based 
on an idea implicitly contained in \cit!{25}, it
was fully proved in \cit!{5} for compact $M$ 
and for general $M$ in \cit!{18}. The clearest 
presentation is in \cit!{19}, section 13. 
If we take a Hilbert space $H$ instead of $N$, then  $B(M,H)$ is 
the classifying space for $\Diff(M)$ if $M$ is compact, and the 
classifying bundle $\Emb(M,H)$ carries also a universal 
connection. This is shown in \cit!{21}. 
 
\subhead\nmb.{2.2} \endsubhead 
If $(N,g)$ is a Riemannian manifold then on the manifold $\Emb(M,N)$  
there is a naturally induced weak Riemannian metric 
given, for  $s_1,s_2\in\Ga_c(f^*TN)$ and $\ph\in \Emb(M,N)$, by
$$ 
G_\phi(s_1,s_2)=\int_Mg(s_1,s_2)\operatorname{vol}(\phi^*g), \quad
\phi \in \Emb(M,N),  
$$ 
where $\operatorname{vol}(g)$ denotes the volume form on $N$ induced
by the Riemannian metric $g$ and $\operatorname{vol}(\phi^*g)$ the
volume form on $M$ induced by the pull back metric $\phi^*g$. The
covariant derivative and curvature of the Levi-Civita connection
induced by $G$ were investigated in
\cit!{4} if $N=\Bbb R^{\dim M+1}$ (endowed with the standard
inner product) and in \cit!{9} for the general case. We shall not
reproduce the general formulae here 
 
This weak Riemannian  
metric is invariant under the action of the diffeomorphism group  
$\Diff(M)$ by composition from the right and hence it induces a  
Riemannian metric on the base manifold $B(M,N)$. 
 
\subhead\nmb.{2.3}. Example \endsubhead 
Let us consider the special case $M=N=\Bbb R$, that is, the space
$\Emb(\Bbb R,\Bbb R)$  of all embeddings of the real line into
itself, which contains the   diffeomorphism group $\Diff(\Bbb R)$ as
an open subset.  The case $M=N=S^1$ is treated in a similar fashion
and the results of this paper are also valid in this situation, where 
$\Emb(S^1,S^1)=\Diff(S^1)$. 
For our purposes, we may restrict attention to
the space of orientation-preserving embeddings, denoted by
$\Emb^{+}(\Bbb R,\Bbb R)$.
The weak Riemannian metric has thus the expression 
$$ 
G_f(h,k) = \int_\Bbb R h(x)k(x)|f'(x)|\,dx,\quad  
f\in\Emb(\Bbb R,\Bbb R), \quad h,k\in C^\infty_c(\Bbb R,\Bbb R). 
$$ 
We shall compute the geodesic equation for this metric by variational  
calculus. The energy of a curve $f$ of embeddings is
$$ 
E(f) = \tfrac12\int_a^b G_f(f_t,f_t) dt  
     = \tfrac12\int_a^b \!\!\int_{\Bbb R} f_t^2f_x\, dx dt. 
$$ 
If we assume that $f(x,t,s)$ is a smooth function and that the
variations are with fixed endpoints, then the derivative with
respect to $s$ of the energy is  
$$\align 
\left. \dd s \right|_0 E(f(\quad,s)) 
&= \left. \dd s \right |_0 \tfrac12\int_a^b \!\!\int_{\Bbb R}
         f_t^2f_x\, dx dt\\  
&= \tfrac12\int_a^b \!\!\int_{\Bbb R}
      (2f_tf_{ts}f_x+f_t^2f_{xs})dx dt\\  
&=  -\tfrac12\int_a^b \!\!\int_{\Bbb R}  
       (2f_{tt}f_{s}f_x+2f_{t}f_{s}f_{tx}+2f_tf_{tx}f_s)dx dt\\ 
&= -\int_a^b \!\!\int_{\Bbb R}  
     \left(f_{tt}+2\frac{f_{t}f_{tx}}{f_x}\right)f_{s}f_xdx dt,\\ 
\endalign$$ 
so that the geodesic equation with its initial data is:  
$$\align 
f_{tt}&=-2\frac{f_{t}f_{tx}}{f_x},\quad  
     f(\quad,0)\in \Emb^+(\Bbb R,\Bbb R), \quad 
     f_t(\quad,0)\in C^\infty_c(\Bbb R,\Bbb R)\tag1\\ 
&=: \Ga_f(f_t,f_t), 
\endalign$$ 
where the Christoffel symbol  
$\Ga:\Emb(\Bbb R,\Bbb R)\x C^\infty_c(\Bbb R,\Bbb R)\x  
C^\infty_c(\Bbb R,\Bbb R)\to C^\infty_c(\Bbb R,\Bbb R)$ 
is given by symmetrisation:
$$ 
\Ga_f(h,k) := -\frac{hk_x+h_xk}{f_x}= -\frac{(hk)_x}{f_x}\,.\tag2 
$$ 
For vector fields $X,Y$ on $\Emb(\Bbb R,\Bbb R)$ the  
covariant derivative is given by the expression  
$\nabla^{\Emb}_XY= dY(X)-\Ga(X,Y)$. 
The Riemannian curvature  
$R(X,Y)Z = (\nabla_X\nabla_Y - \nabla_Y\nabla_X -\nabla_{[X,Y]})Z$ 
is then determined in terms of the Christoffel form by
$$\align 
R_f(h,k)\ell &= -d\Ga(f)(h)(k,\ell) + d\Ga(f)(k)(h,\ell) + 
     \Ga_f(h,\Ga_f(k,\ell)) - \Ga_f(k,\Ga_f(h,\ell))\\ 
&= -\frac{h_x(k\ell)_x}{f_x^2} +\frac{k_x(h\ell)_x}{f_x^2} 
     +\frac{\left( h\frac{(k\ell)_x}{f_x} \right)_x}{f_x} 
     -\frac{\left( k\frac{(h\ell)_x}{f_x} \right)_x}{f_x}\\ 
&= \frac1{f_x^3}\Bigl( 
     f_{xx}h_xk\ell - f_{xx}hk_x\ell  
     +f_xhk_{xx}\ell - f_xh_{xx}k\ell 
     +2f_xhk_x\ell_x -2f_xh_xk\ell_x\Bigr)\tag3 
\endalign$$ 
The geodesic equation can be solved in the following way: 
if instead of the obvious framing we change variables to  
$T\Emb= \Emb\x C^\infty_c\ni (f,h) \mapsto (f,hf_x^2)=:(f,F)$ 
then the geodesic equation becomes  
$F_t=\dd t(f_tf_x^2)=f_x^2(f_{tt}+2\frac{f_tf_{tx}}{f_x})=0,$ 
so that $F=f_tf_x^2$ is constant in $t$, or 
$f_t(x,t)f_x(x,t)^2=f_t(x,0)f_x(x,0)^2$. From
here, a standard separation of variables
argument gives the solution; it blows up in
finite time for most initial conditions.  

Now let us consider the trivialisation of $T\Emb(\Bbb R,\Bbb R)$ by  
right translation (this is most useful for $\Diff(\Bbb R)$). The
derivative of the inversion $\operatorname{Inv}:g\mapsto  g\i$ is given by
$$
T_g(\operatorname{Inv})h = - T(g\i)\o h \o g\i = 
\frac{h\o g\i}{g_x\o g\i}\quad\text{ for }\quad g \in
\Emb(\Bbb R, \Bbb R),\, h \in C^\infty_c(\Bbb R,\Bbb R).
$$
Defining  
$$
u:= f_t\o f\i,\quad\text{ or, in more detail, }\quad
u(x,t)=f_t(f(\quad,t)\i(x), t) ,
$$
we have
$$\align 
u_x &= (f_t \o f\i)_x = (f_{tx}\o f\i)\frac1{f_x\o f\i}
     = \frac{f_{tx}}{f_x} \o f\i,\\ 
u_t &= (f_t \o f\i)_t = f_{tt}\o f\i + (f_{tx}\o f\i)(f\i)_t \\
&= f_{tt}\o f\i + (f_{tx}\o f\i)\frac{1}{f_x\ f\i}(f_t\ f\i)
\endalign$$
which, by \thetag1 and the first equation becomes
$$
u_t = f_{tt}\o f\i - \left( \frac{f_{tx}f_t}{f_x} \right)\o f\i
     = -3\left(\frac{f_{tx}f_t}{f_x}\right)\o f\i 
     = -3 u_xu.
$$
The geodesic equation on $\Emb(\Bbb R, \Bbb R)$ in right
trivialization, that is, in Eulerian formulation, is hence  
$$ 
u_t = -3 u_xu \,,\tag 4 
$$  
which is just Burgers' equation. 
 
\head\totoc\nmb0{3}. Right invariant Riemannian metrics on Lie groups  
\endhead 
 
\subhead\nmb.{3.1}. Geodesics of a right invariant metric on a Lie  
group \endsubhead 
Let $G$ be a Lie group which may be infinite dimensional, with Lie  
algebra $\g$. 
Let $\mu:G\x G\to G$ be the multiplication, let $\mu_x$ be left  
translation and $\mu^y$ be right translation, 
given by $\mu_x(y)=\mu^y(x)=xy=\mu(x,y)$. We also need the right  
Maurer-Cartan form $\ka=\ka^r\in\Om^1(G,\g)$, given by  
$\ka_x(\xi):=T_x(\mu^{x\i})\cdot \xi$. It satisfies the right
Maurer-Cartan equation $d\ka-\tfrac12[\ka,\ka]_\wedge=0$, where
$[\;,\;]_\wedge$ denotes the wedge product of $\g$-valued forms on
$G$ induced by the Lie bracket. Note that
$\tfrac12[\ka,\ka]_\wedge (\xi,\et) = [\ka(\xi),\ka(\et)]$.
 
Let $\langle \;,\;\rangle:\g\x\g\to\Bbb R$ be a positive  
definite bounded (weak) inner product. Then  
$$ 
G_x(\xi,\et)=\langle T(\mu^{x\i})\cdot\xi,
     T(\mu^{x\i})\cdot\et)\rangle =  
     \langle \ka(\xi),\ka(\et)\rangle  \tag1
$$ 
is a right invariant (weak) Riemannian metric on $G$, and any
(weak) right invariant bounded Riemannian metric is of this form, for
suitable  $\langle \;,\;\rangle$. 
 
Let $g:[a,b]\to G$ be a smooth curve.  
The velocity field of $g$, viewed in the right trivializations, 
coincides  with the right logarithmic derivative
$T(\mu^{g\i})\cdot \partial_t g =  
\ka(\partial_t g) = (g^*\ka)(\partial_t)$, where 
$\partial_t=\frac{\partial}{\partial t}$. 
The energy of the curve $g(t)$ is given by  
$$ 
E(g) = \tfrac12\int_a^bG_g(g',g')dt = \tfrac12\int_a^b 
     \langle (g^*\ka)(\partial_t),(g^*\ka)(\partial_t)\rangle\, dt. 
$$ 
For a variation $g(t,s)$ with fixed endpoint we have then, using the  
right Maurer-Cartan equation and integration by parts, 
$$\align 
\partial_sE(g) &= \tfrac12\int_a^b2 
     \langle \partial_s(g^*\ka)(\partial_t),\,
                            (g^*\ka)(\partial_t)\rangle\, dt\\ 
&= \int_a^b \langle \partial_t(g^*\ka)(\partial_s) - 
         d(g^*\ka)(\partial_t,\partial_s),\,
                (g^*\ka)(\partial_t)\rangle\,dt\\   
&= \int_a^b \left(-\langle (g^*\ka)(\partial_s),\, 
     \partial_t(g^*\ka)(\partial_t)\rangle - 
     \langle [(g^*\ka)(\partial_t),(g^*\ka)(\partial_s)],\, 
     (g^*\ka)(\partial_t)\rangle\right)\, dt\\ 
&= -\int_a^b \langle (g^*\ka)(\partial_s),\, 
     \partial_t(g^*\ka)(\partial_t) + 
     \ad((g^*\ka)(\partial_t))^{\top}
           ((g^*\ka)(\partial_t))\rangle\, dt\\ 
\endalign$$ 
where $\ad((g^*\ka)(\partial_t))^{\top}:\g\to\g$ is the adjoint of  
$\ad((g^*\ka)(\partial_t))$ with respect to the inner product  
$\langle \;,\; \rangle$. In infinite dimensions one also has to  
check the existence of this adjoint. 
In terms of the right logarithmic derivative $u:[a,b]\to \g$ of  
$g:[a,b]\to G$, given by  
$u(t):= g^*\ka(\partial_t) = T_{g(t)}(\mu^{g(t)\i})\cdot g'(t)$, 
the geodesic equation has the expression 
$$ 
u_t = - \ad(u)^{\top}u\,. \tag2
$$
This is, of course, just the Euler-Poincar\'e equation for right
invariant systems using the Lagrangian given by the kinetic energy
(see \cit!{15}, section 13) and the above derivation is done directly
without invoking this theorem. 
 
\subhead\nmb.{3.2}. The covariant derivative \endsubhead 
Our next aim is to derive the Riemannian curvature and for that we  
develop the basis-free version of Cartan's method of moving frames 
in this setting, which also works in infinite dimensions. 
The right trivialization, or framing, $(\ka,\pi_G):TG\to \g\x G$  
induces the isomorphism $R:C^\infty(G,\g)\to \X(G)$, given by  
$R(X)(x):= R_X(x):=T_e(\mu^x)\cdot X(x)$, for $X\in C^\infty(G,\g)$ and
$x\in G$. Here $\X(G):=\Ga(TG)$ denote the Lie algebra of all vector 
fields. For the Lie bracket and the Riemannian metric we have 
$$\gather 
[R_X,R_Y] = R(-[X,Y]_\g + dY\cdot R_X - dX\cdot R_Y),\tag1\\ 
R\i[R_X,R_Y] = -[X,Y]_\g + R_X(Y) - R_Y(X),\\ 
G_x(R_X(x),R_Y(x)) = \langle X(x),Y(x)\rangle\,,\, x\in G. 
\endgather$$ 
In the sequel we shall compute in $C^\infty(G,\g)$ instead of 
$\X(G)$. In particular, we shall use the convention
$$
\nabla_XY := R\i(\nabla_{R_X}R_Y)\quad\text{ for }X,Y\in C^\infty(G,\g).
$$
to express the Levi-Civita covariant derivative. 
 
\proclaim{Lemma} 
Assume that for all $\xi\in\g$ the adjoint $\ad(\xi)^\top$  
with respect to the inner  
product $\langle \;,\;\rangle$ exists and that  
$\xi\mapsto \ad(\xi)^\top$ is bounded. 
Then the Levi-Civita covariant derivative of the metric  
2.1\thetag1 exists and is given for any $X,Y \in C^\infty(G,\g)$  in
terms of the isomorphism $R$ by
$$
\nabla_XY= dY.R_X + \tfrac12\ad(X)^\top Y 
     +\tfrac12\ad(Y)^\top X - \tfrac12\ad(X)Y.
\tag2$$
\endproclaim 
 
\demo{Proof}  
Easy computations show that this formula satisfies the axioms of a
covariant derivative, that relative to it the Riemannian
metric is covariantly constant, since 
$$
R_X\langle Y,Z\rangle = \langle dY.R_X,Z\rangle + \langle Y,dZ.R_X\rangle
     = \langle \nabla_XY,Z\rangle + \langle Y,\nabla_XZ\rangle,
$$
and that it is torsion free, since 
$$
\nabla_XY-\nabla_YX + [X,Y]_\g - dY.R_X + dX.R_Y = 0.\qed
$$
\enddemo
\smallskip 
 
For $\xi \in \g$ define $\al(\xi):\g\to\g$ by
$\al(\xi)\eta:=\ad(\eta)^\top \xi$. With this notation, the previous
lemma states that for all $X\in C^\infty(G,\g)$ the covariant
derivative of the Levi-Civita connection has the expression
$$
\nabla_X = R_X + \tfrac12\ad(X)^\top +\tfrac12\al(X) - 
\tfrac12\ad(X).\tag3
$$
 
\subhead\nmb.{3.3}. The curvature \endsubhead 
First note that we have the following relations: 
$$\alignat2 
[R_X,\ad(Y)] &= \ad(R_X(Y)),&\quad  
[R_X,\al(Y)] &= \al(R_X(Y)),\tag1\\ 
[R_X,\ad(Y)^\top] &= \ad(R_X(Y))^\top,&\quad 
[\ad(X)^\top,\ad(Y)^\top] &= -\ad([X,Y]_\g)^\top.  
\endalignat$$ 
The Riemannian curvature is then computed by  
$$\align 
&\Cal R(X,Y) =  
[\nabla_X,\nabla_Y]-\nabla_{-[X,Y]_\g+R_X(Y)-R_Y(X)}\tag2\\ 
&= [R_X+\tfrac12\ad(X)^\top+\tfrac12\al(X)-\tfrac12\ad(X), 
     R_Y+\tfrac12\ad(Y)^\top+\tfrac12\al(Y)-\tfrac12\ad(Y)]\\ 
&\quad - R_{-[X,Y]_\g + R_X(Y) - R_Y(X)} 
     -\tfrac12\ad(-[X,Y]_\g + R_X(Y) - R_Y(X))^\top\\ 
&\quad -\tfrac12\al(-[X,Y]_\g + R_X(Y) - R_Y(X)) 
     +\tfrac12\ad(-[X,Y]_\g + R_X(Y) - R_Y(X)) \\ 
&= -\tfrac14[\ad(X)^\top+\ad(X),\ad(Y)^\top+\ad(Y)]\\ 
&\quad+\tfrac14[\ad(X)^\top-\ad(X),\al(Y)]  
      +\tfrac14[\al(X),\ad(Y)^\top-\ad(Y)] \\ 
&\quad+\tfrac14[\al(X),\al(Y)] +\tfrac12\al([X,Y]_\g).
\endalign$$ 
If we plug in all definitions and use 4 times 
the Jacobi identity we  
get the following expression 
$$\align 
&\langle 4\Cal R(X,Y)Z,U\rangle = 
+ 2\langle [X,Y],[Z,U] \rangle                   
- \langle [Y,Z],[X,U] \rangle                   
+ \langle [X,Z],[Y,U] \rangle\\ 
&- \langle Z,[U,[X,Y]] \rangle  
+ \langle U,[Z,[X,Y]] \rangle  
- \langle Y,[X,[U,Z]] \rangle                   
- \langle X,[Y,[Z,U]] \rangle\\ 
&+ \langle \ad(X)^\top Z,\ad(Y)^\top U \rangle   
+ \langle \ad(X)^\top Z,\ad(U)^\top Y \rangle   
+ \langle \ad(Z)^\top X,\ad(Y)^\top U \rangle\\   
&- \langle \ad(U)^\top X,\ad(Y)^\top Z \rangle 
- \langle \ad(Y)^\top Z,\ad(X)^\top U \rangle 
- \langle \ad(Z)^\top Y,\ad(X)^\top U \rangle\\ 
&- \langle \ad(U)^\top X,\ad(Z)^\top Y \rangle   
+ \langle \ad(U)^\top Y,\ad(Z)^\top X \rangle. 
\endalign$$ 
 
\subhead\nmb.{3.4}. Jacobi fields, I \endsubhead 
We compute first the Jacobi equation directly via variations of 
geodesics.  
So let $g:\Bbb R^2\to G$ be smooth, $t\mapsto g(t,s)$ a geodesic for  
each $s$. Let again $u=\ka(\partial_t g) = (g^*\ka)(\partial_t)$ be  
the velocity field along the geodesic in right trivialization which  
satisfies the geodesic equation  
$u_t=-\ad(u)^\top u$. 
Then $y:= \ka(\partial_s g)=(g^*\ka)(\partial_s)$ is the Jacobi field  
corresponding to this variation, written in the right trivialization. 
 From the right Maurer-Cartan equation we then have: 
$$\align 
y_t &= \partial_t (g^*\ka)(\partial_s)  
     = d(g^*\ka)(\partial_t,\partial_s)  
     + \partial_s (g^*\ka)(\partial_t) + 0\\ 
&= [(g^*\ka)(\partial_t),(g^*\ka)(\partial_s)]_{\frak g} + u_s\\ 
&= [u,y] +u_s. 
\endalign$$ 
Using the geodesic equation, the definition of $\al$, and
the fourth relation in \nmb!{3.3}.\thetag1, this identity implies
$$\align 
u_{st} &= u_{ts} = \partial_s u_t = -\partial_s (\ad(u)^\top u)  
     = -\ad(u_s)^\top u - \ad(u)^\top u_s\\ 
&= -\ad(y_t+[y,u])^\top u - \ad(u)^\top (y_t+[y,u])\\ 
&= -\al(u)y_t -\ad([y,u])^\top u  
     - \ad(u)^\top y_t - \ad(u)^\top ([y,u])\\ 
&= - \ad(u)^\top y_t -\al(u)y_t +[\ad(y)^\top,\ad(u)^\top]u  
     - \ad(u)^\top\ad(y)u\,.\\ 
\endalign$$ 
Finally we get the Jacobi equation as 
$$\align 
y_{tt} &= [u_t,y] + [u,y_t] + u_{st}\\ 
&= \ad(y)\ad(u)^\top u + \ad(u)y_t - \ad(u)^\top y_t \\ 
&\qquad-\al(u)y_t +[\ad(y)^\top,\ad(u)^\top]u - 
\ad(u)^\top\ad(y)u\,,\\ 
y_{tt}&= [\ad(y)^\top+\ad(y),\ad(u)^\top]u   
     - \ad(u)^\top y_t -\al(u)y_t + \ad(u)y_t\,.\tag1 
\endalign$$ 
 
\subhead\nmb.{3.5}. Jacobi fields, II \endsubhead 
Let $y$ be a Jacobi field along a geodesic 
$g$ with right trivialized  
velocity field $u$. Then $y$ should satisfy the analogue of the 
finite dimensional Jacobi equation 
$$ 
\nabla_{\partial_t}\nabla_{\partial_t}y + \Cal R(y,u)u = 0 
$$ 
We want to show that this leads to same equation 
as \nmb!{3.4}.\thetag1. 
First note that from \nmb!{3.2}.\thetag2 we have 
$$ 
\nabla_{\partial_t}y  
     = y_t +\tfrac12\ad(u)^\top y + \tfrac12\al(u)y - 
\tfrac12\ad(u)y 
$$ 
so that, using $u_t=-\ad(u)^\top u$, we get: 
$$\align 
\nabla_{\partial_t}\nabla_{\partial_t}y  
&= \nabla_{\partial_t}\Bigl(y_t +\tfrac12\ad(u)^\top y + 
\tfrac12\al(u)y -  
\tfrac12\ad(u)y\Bigr)\\ 
\allowdisplaybreak
&= y_{tt} +\tfrac12\ad(u_t)^\top y +\tfrac12\ad(u)^\top y_t  
     + \tfrac12\al(u_t)y \\ 
&\quad+ \tfrac12\al(u)y_t - \tfrac12\ad(u_t)y  
     -\tfrac12\ad(u)y_t\\ 
&\quad+\tfrac12\ad(u)^\top\Bigl(y_t+\tfrac12\ad(u)^\top y  
     + \tfrac12\al(u)y - \tfrac12\ad(u)y\Bigr) \\ 
&\quad+\tfrac12\al(u)\Bigl(y_t+\tfrac12\ad(u)^\top y 
     + \tfrac12\al(u)y - \tfrac12\ad(u)y\Bigr) \\ 
&\quad-\tfrac12\ad(u)\Bigl(y_t+\tfrac12\ad(u)^\top y 
     + \tfrac12\al(u)y - \tfrac12\ad(u)y\Bigr) \\ 
\allowdisplaybreak
&= y_{tt}  + \ad(u)^\top y_t + \al(u)y_t - \ad(u)y_t\\ 
&\quad-\tfrac12\al(y)\ad(u)^\top u  
      - \tfrac12\ad(y)^\top\ad(u)^\top u  
      - \tfrac12\ad(y)\ad(u)^\top u \\ 
&\quad+\tfrac12\ad(u)^\top\Bigl(\tfrac12\al(y)u  
     + \tfrac12\ad(y)^\top u + \tfrac12\ad(y)u\Bigr) \\ 
&\quad+\tfrac12\al(u)\Bigl(\tfrac12\al(y)u 
     + \tfrac12\ad(y)^\top u + \tfrac12\ad(y)u\Bigr) \\ 
&\quad-\tfrac12\ad(u)\Bigl(\tfrac12\al(y)u 
     + \tfrac12\ad(y)^\top u + \tfrac12\ad(y)u\Bigr)\,. \\ 
\endalign$$ 
In the second line of the last expression we use  
$$ 
-\tfrac12\al(y)\ad(u)^\top u= -\tfrac14\al(y)\ad(u)^\top u  
     -\tfrac14\al(y)\al(u)u 
$$ 
and similar forms for the other two terms to get: 
$$\align 
\nabla_{\partial_t}\nabla_{\partial_t}y  
&= y_{tt}  + \ad(u)^\top y_t + \al(u)y_t - \ad(u)y_t\\ 
&\quad+\tfrac14[\ad(u)^\top,\al(y)]u  
     + \tfrac14[\ad(u)^\top,\ad(y)^\top]u  
     + \tfrac14[\ad(u)^\top,\ad(y)]u \\ 
&\quad+\tfrac14[\al(u),\al(y)]u 
     + \tfrac14[\al(u),\ad(y)^\top]u  
     + \tfrac14[\al(u),\ad(y)]u \\ 
&\quad-\tfrac14[\ad(u),\al(y)]u 
     - \tfrac14[\ad(u),\ad(y)^\top + \ad(y)]u, \\ 
\endalign$$ 
where in the last line we also used $\ad(u)u=0$. 
We now compute the curvature term using \nmb!{3.3}.\thetag2: 
$$\align 
\Cal R(y,u)u 
&= -\tfrac14[\ad(y)^\top+\ad(y),\ad(u)^\top+\ad(u)]u\\ 
&\quad+\tfrac14[\ad(y)^\top-\ad(y),\al(u)]u  
      +\tfrac14[\al(y),\ad(u)^\top-\ad(u)]u \\ 
&\quad+\tfrac14[\al(y),\al(u)] +\tfrac12\al([y,u])u\\ 
&= -\tfrac14[\ad(y)^\top+\ad(y),\ad(u)^\top]u 
      -\tfrac14[\ad(y)^\top+\ad(y),\ad(u)]u\\ 
&\quad+\tfrac14[\ad(y)^\top,\al(u)]u  
      -\tfrac14[\ad(y),\al(u)]u  
      +\tfrac14[\al(y),\ad(u)^\top-\ad(u)]u \\ 
&\quad+\tfrac14[\al(y),\al(u)]u +\tfrac12\ad(u)^\top\ad(y)u\,.\\ 
\endalign$$ 
Summing up we get 
$$\align 
\nabla_{\partial_t}\nabla_{\partial_t}y + \Cal R(y,u)u  
&= y_{tt}  + \ad(u)^\top y_t + \al(u)y_t - \ad(u)y_t \\ 
&\quad-\tfrac12[\ad(y)^\top+\ad(y),\ad(u)^\top]u \\ 
&\quad+ \tfrac12[\al(u),\ad(y)]u +\tfrac12\ad(u)^\top\ad(y)u\,. 
\endalign$$ 
Finally we need the following computation using \nmb!{3.3}.\thetag1: 
$$\align 
\tfrac12[\al(u),\ad(y)]u  
&= \tfrac12\al(u)[y,u] - \tfrac12\ad(y)\al(u)u\\ 
&= \tfrac12\ad([y,u])^\top u - \tfrac12\ad(y)\ad(u)^\top u\\ 
&= -\tfrac12[\ad(y)^\top,\ad(u)^\top] u -
       \tfrac12\ad(y)\ad(u)^\top u\,. 
\endalign$$ 
Inserting we get the desired result: 
$$\align 
\nabla_{\partial_t}\nabla_{\partial_t}y + \Cal R(y,u)u  
&= y_{tt}  + \ad(u)^\top y_t + \al(u)y_t - \ad(u)y_t \\ 
&\quad-[\ad(y)^\top+\ad(y),\ad(u)^\top]u. 
\endalign$$ 

\subhead\nmb.{3.6}. The weak symplectic structure on 
the space of Jacobi  
fields \endsubhead 
Let us assume now that the geodesic equation in $\g$ 
$$ 
u_t=-\ad(u)^\top u  
$$ 
admits a unique solution for some time 
interval, depending smoothly on the  
choice of the initial value $u(0)$. 
Furthermore we assume that $G$ is  
a regular Lie group (see \cit!{13}, 5.3) 
so that each smooth curve $u$ in $\g$  
is the right logarithmic derivative of a smooth curve  
$g$ in $G$ which depends smoothly on $u$, 
so that $u=(g^*\ka)(\partial_t)$.
Furthermore we have to assume 
that the Jacobi equation along $u$ 
admits a unique solution for some time,
depending smoothly on the initial values $y(0)$ and $y_t(0)$.
These are non-trivial assumptions: 
in \cit!{13}, 2.4 there are examples of
ordinary linear differential equations 
`with constant coefficients' which
violate existence or uniqueness. 
These assumptions have to be 
checked in the special situations.
Then the space $\Cal J_u$ of all Jacobi fields 
along the geodesic $g$ described by $u$
is isomorphic to the space $\g\x \g$ of all initial data.
 
There is the well known symplectic structure on the space $\Cal J_u$  
of all Jacobi fields along a fixed geodesic with velocity field $u$, see
e.g\. \cit!{11}, II, p.70.  
It is given by the following expression which is constant in time $t$: 
$$\align 
\si(y,z) :&= \langle y,\nabla_{\partial_t}z \rangle  
     - \langle \nabla_{\partial_t}y,z \rangle \\ 
&= \langle y, z_t + \tfrac12\ad(u)^\top z 
          + \tfrac12\al(u)z - \tfrac12\ad(u)z \rangle \\ 
&\qquad-\langle y_t + \tfrac12\ad(u)^\top y  
          + \tfrac12\al(u)y - \tfrac12\ad(u)y, z  \rangle\\ 
&= \langle y,z_t \rangle - \langle y_t, z \rangle  
     + \langle [u,y],z \rangle - \langle y,[u,z] \rangle  
     - \langle [y,z],u \rangle \\ 
&= \langle y,z_t-\ad(u)z+\tfrac12\al(u)z \rangle  
     -\langle y_t-\ad(u)y+\tfrac12\al(u)y, z \rangle.  
\endalign$$ 
It is worth while to check directly from the Jacobi field equation 
\nmb!{3.4}.\thetag1 that $\si(y,z)$ is indeed constant in $t$. 
Clearly $\si$ is a weak symplectic structure on the relevant
vector space $\Cal J_u\cong \g\x \g$, i.e., $\si$ gives an
injective (but in general not surjective) linear
mapping $\Cal J_u\to \Cal J_u^*$. This is seen most
easily by writing 
$$
\si(y,z)=\langle y, z_t - \Ga_g(u,z)\rangle|_{t=0} 
	- \langle y_t-\Ga_g(u,y),z\rangle|_{t=0}
$$
which is induced from the standard symplectic structure on 
$\g\x \g^*$ by applying first the automorphism 
$(a,b)\mapsto (a,b-\Ga_g(u,a))$  to $\g\x\g$ and then by injecting 
the second factor $\g$ into its dual $\g^*$.

For regular (infinite dimensional) Lie groups variations of geodesics
exist, but there is no general theorem stating that they are uniquely
determined by $y(0)$ and $y_t(0)$. For concrete regular Lie groups, this
needs to be shown directly.   
 
\head\totoc\nmb0{4}. The diffeomorphism group of the circle revisited  
\endhead 
 
\subhead\nmb.{4.1}. Geodesics and curvature \endsubhead 
We consider again the Lie groups $\Diff(\Bbb R)$ and  
$\Diff(S^1)$ with Lie algebras $\X_c(\Bbb R)$ and $\X(S^1)$ where the  
Lie bracket $[X,Y]=X'Y-XY'$ is the negative of the usual one. For the  
inner product $\langle X,Y\rangle=\int X(x)Y(x)\,dx$ 
integration by parts gives 
$$ 
\langle [X,Y],Z\rangle = \int_\Bbb R (X'YZ-XY'Z)dx  
     = \int_\Bbb R (2X'YZ+XYZ')dx = \langle Y,\ad(X)^\top Z\rangle,  
$$ 
which in turn gives rise to  
$$\align 
\ad(X)^\top Z &= 2X'Z+XZ', \tag1\\  
\al(X)Z &= 2Z'X+ZX' ,\tag2\\ 
(\ad(X)^\top+\ad(X)) Z &= 3X'Z ,\tag3\\ 
(\ad(X)^\top-\ad(X)) Z &= X'Z + 2XZ' = \al(X)Z . \tag4  
\endalign$$ 
Equation \thetag4 states that $-\frac12\al(X)$ 
is the skew-symmetrization of  
$\ad(X)$ with respect to to the inner product  
$\langle\quad,\quad\rangle$. From the 
theory of symmetric spaces one  
then expects that $-\frac12\al$ is a 
Lie algebra homomorphism and  
indeed one can check that  
$$ 
-\tfrac12\al([X,Y]) = \left[-\tfrac12\al(X),
                       -\tfrac12\al(Y)\right] 
$$ 
holds for any vector fields $X, Y$. From 
\thetag1 we get the same geodesic equation as in  
\nmb!{2.3}\thetag4, namely Burgers' equation: 
$$ 
u_t = -\ad(u)^\top u = -3u_xu. 
$$ 
Using the above relations and the general curvature formula  
\nmb!{3.3}.\thetag2, we get
$$\align 
\Cal R(X,Y)Z &= -X''YZ +XY''Z -2X'YZ' +2XY'Z'  
             = -2[X,Y]Z' - [X,Y]'Z \\ 
&= -\al([X,Y])Z. \tag5
\endalign$$ 
If we change the framing of the tangent bundle: 
$$ 
X=h\o f\i, \quad X' = \left(\frac{h_x}{f_x}\right)\o f\i,  
\quad X'' = \left(\frac{h_{xx}f_x-h_xf_{xx}}{f_x^3}\right)\o f\i,  
$$  
and similarly for $Y=k\o f\i$ and $Z=\ell\o f\i$, for
$h, k, \ell \in  C^\infty_c(\Bbb R,\Bbb R)$ or 
$C^\infty(S^1,\Bbb R)$, then  
$(\Cal R(X,Y)Z)\o f$ given by \thetag5 coincides with formula
\nmb!{2.3}.\thetag3 for the curvature.  
 
\subhead\nmb.{4.2}. Jacobi fields \endsubhead 
A Jacobi field $y$ along a geodesic $g$ with velocity field $u$  
is a solution of the  
partial differential equation \nmb!{3.4}.\thetag1, which in our case  
becomes: 
$$\align 
y_{tt}&= [\ad(y)^\top+\ad(y),\ad(u)^\top]u  
     - \ad(u)^\top y_t -\al(u)y_t + \ad(u)y_t \tag1\\ 
&= - 3u^2y_{xx} -4 uy_{tx} -2u_xy_t\\ 
u_t &= -3 u_xu. 
\endalign$$ 
Since the geodesic equation has solutions, locally in time (see the  
argument in \nmb!{2.3}) it is to be expected that 
the space of all Jacobi fields exists and is isomorphic to the space 
of all initial data $(y(0),y_t(0))\in C^\infty(S^1,\Bbb R)^2$ or 
$C^\infty_c(\Bbb R,\Bbb R)^2$,   
respectively. 
The weak symplectic structure on it is given by \nmb!{3.6}: 
$$\align 
\si(y,z) &= \langle y, z_t-\tfrac12u_xz+2uz_x \rangle 
     - \langle y_t-\tfrac12u_xy+2uy_x, z \rangle \\ 
&= \int_{S^1\text{ \!\!or }\Bbb R} (yz_t-y_tz+2u(yz_x-y_xz))\,dx. \tag2 
\endalign$$ 
 
\head\totoc\nmb0{5}. The Virasoro-Bott group and the Korteweg-de  
Vries-equation \endhead 
 
\subhead\nmb.{5.1}. Geodesics on the Virasoro-Bott group  \endsubhead 
For $\ph\in\Diff^+(S^1)$ let $\ph':S^1\to \Bbb R^+$ be the mapping 
given by $T_x\ph\cdot\partial_x=\ph'(x)\partial_x$. 
Then  
$$\gather 
c:\Diff^+(S^1)\x\Diff^+(S^1)\to \Bbb R\\ 
c(\ph,\ps):=\int_{S^1}\log(\ph\o\ps)'d\log\ps'  
     = \int_{S^1}\log(\ph'\o\ps)d\log\ps' 
\endgather$$ 
satisfies $c(\ph,\ph\i)=0$ and is a 
smooth group cocycle, called the Bott cocycle.
The corresponding central extension group
$S^1\x_c\Diff^+(S^1)$, called the Virasoro-Bott group, is
a trivial $S^1$-bundle  
$S^1\x\Diff^+(S^1)$ that becomes a regular Lie  
group relative to the operations  
$$\binom{\ph}{\al}\binom{\ps}{\be}  
     =\binom{\ph\o\ps}{\al\be\,e^{2\pi ic(\ph,\ps)}},\quad 
\binom{\ph}{\al}\i=\binom{\ph\i}{\al\i}\quad 
\ph, \ps \in \Diff^+(S^1),\; \al,\be \in S^1 . 
$$ 
The Lie algebra of this Lie group is the central extension 
$\Bbb R\x_\om \X(S^1)$ of $\X(S^1)$ induced by the Gelfand-Fuchs
Lie algebra cocycle $\om:\X(S^1)\x \X(S^1)\to \Bbb R$ 
$$ 
\om(h,k)=\om(h)k=\int_{S^1}h'dk'=\int_{S^1}h'k''dx =  
\tfrac12\int_{S^1}\det\pmatrix h'& k'\\ h''&k''\endpmatrix\,dx, 
$$   
a generator of the 1-dimensional bounded Chevalley cohomology  
$H^2(\X(S^1),\Bbb R)$. Thus the bracket on $\Bbb R\x_\om \X(S^1)$
is given by  
$$ 
\left[\binom{h}{a},\binom{k}{b}\right] = \binom{h'k-hk'}{\om(h,k)},
\quad h, k \in \X(S^1),\; a,b \in \Bbb R. 
$$ 
Note that the Lie algebra cocycle $\om$ makes sense  
on the Lie algebra $\X_c(\Bbb R)$ of all vector fields with compact  
support on $\Bbb R$, but that it does not integrate to 
a group cocycle on  
$\Diff(\Bbb R)$. The subsequent considerations also make sense on  
$\X_c(\Bbb R)$. Recall also that $H^2(\X_c(M),\Bbb R)=0$ for each
finite dimensional manifold of dimension $\ge 2$ (see \cit!{7}),  
which blocks the way to  
find a higher dimensional analog of the Korteweg -- de Vries  
equation in a way similar to that sketched below. 
 
We shall use the $L^2$-inner product on $\Bbb R\x_\om \X(S^1)$: 
$$ 
\left\langle \binom{h}{a},\binom{k}{b}\right\rangle :=  
\int_{S^1}hk\,dx + ab. 
$$ 
Integrating by parts we get 
$$\align 
\left\langle \ad\binom{h}{a}\binom{k}{b},\binom{\ell}{c}\right\rangle  
&= \left\langle \binom{h'k-hk'}{\om(h,k)}, 
     \binom{\ell}{c}\right\rangle\\ 
&= \int_{S^1}(h'k\ell-hk'\ell+ch'k'')\,dx 
= \int_{S^1}(2h'\ell+h\ell'+ch''')k\,dx\\ 
&= \left\langle \binom{k}{b}, 
     \ad\binom{h}{a}^\top\binom{\ell}{c}\right\rangle, 
          \quad {\text {where}}\\ 
\ad\binom{h}{a}^\top\binom{\ell}{c}&=\binom{2h'\ell+h\ell'+ch'''}{0}. 
\endalign$$ 
Using matrix notation we get therefore (where $\partial :=
\partial_x$) 
$$\align 
\ad\binom{h}{a} &= \pmatrix h' - h\partial & 0 \\  
                                   \om(h) & 0 \endpmatrix \\ 
\ad\binom{h}{a}^\top &= \pmatrix 2h'+h\partial & h'''\\ 
                                      0 & 0 \endpmatrix  \\ 
\al\binom{h}{a} &=\ad\binom{\quad}{\quad}^\top\binom{h}{a}  
     = \pmatrix h'+2h\partial +a\partial^3 & 0 \\ 
                                      0 & 0 \endpmatrix  \\ 
\ad{\binom{h}{a}}^\top + \ad\binom{h}{a} 
     &= \pmatrix 3h' & h''' \\ 
                \om(h) & 0 \endpmatrix  \\ 
\ad\binom{h}{a}^\top - \ad\binom{h}{a} 
     &= \pmatrix h'+2h\partial & h''' \\ 
                -\om(h) & 0 \endpmatrix.  \\ 
\endalign$$ 
 Formula \nmb!{3.1}\thetag2 gives the geodesic equation on the  
Virasoro-Bott group:  
$$ 
\binom{u_t}{a_t}=-\ad\binom{u}{a}^\top\binom{u}{a} 
     =\binom{-3u'u-au'''}{0}. 
$$ 
Thus $a$ is a constant in time and the geodesic equation  
is hence the periodic Korteweg-de Vries equation 
$$ 
u_t+3u_xu+au_{xxx}=0. 
$$ 
Had we worked on $\X_c(\Bbb R)$ we would have obtained the usual 
Korteweg-de Vries   equation. The derivation above is direct and
does not use the Euler-Poincar\'e equations; for a derivation of the
Korteweg-de Vries   equation from this point of view 
see \cit!{15}, section 13.8.
 
\subhead\nmb.{5.2}. The curvature \endsubhead 
The computation of the curvature at the identity element has been done
independently by Misiolek \cit!{22}; our results of course agree. Here
we proceed with a completely general computation that takes advantage
of the formalism introduced so far. Inserting the matrices of
differential-   and integral   operators
$\ad\tbinom{h}{a}^\top$,
$\al\tbinom{h}{a}$, and  
$\ad\tbinom{h}{a}$ etc\. 
given above into formula \nmb!{3.3}\thetag2 and recalling that the matrix  
is applied to vectors of the form $\binom{\ell}{c}$, where $c$ is a  
constant, we see that  
$4\Cal R\left(\tbinom{h_1}{a_1},\tbinom{h_2}{a_2}\right)$ is the following  
$2\x 2$-matrix whose entries are differential- and integral  
operators:  
$$ 
\pmatrix 
{\aligned 
&4(h_1h_2''-h_1''h_2) +2(a_1h_2^{(4)}-a_2h_1^{(4)})\\ 
+(&8(h_1h_2'-h_1'h_2)+10(a_1h_2'''-a_2h_1'''))\partial\\ 
+&18(a_1h_2''-a_2h_1'')\partial^2\\ 
+&(12(a_1h_2'-a_2h_1')+2\om(h_1,h_2))\partial^3\\ 
-&h_1'''\om(h_2)+h_2'''\om(h_1) 
\endaligned} 
&\qquad& 
{\aligned  
&2(h_1'''h_2' - h_1'h_2''')\\ 
+&2(h_1h_2^{(4)}-h_1^{(4)}h_2)\\  
+&(a_1h_2^{(6)}-a_2h_1^{(6)})\endaligned} \\ 
\vphantom{\int_A^A} & & \\ 
{\aligned 
&\om(h_2)(4h_1'+2h_1\partial+a_1\partial^3)\\ 
-&\om(h_1)(4h_2'+2h_2\partial+a_2\partial^3) 
\endaligned} 
&& 0 \endpmatrix .
$$
Therefore, $4\Cal R\left(\tbinom{h_1}{a_1},\tbinom{h_2}{a_2}\right)
\tbinom{h_3}{a_3}$ has the following expression
$$
\pmatrix
{\aligned
&4(h_1h_2'' - h_1''h_2)h_3 + 2(a_1h_2^{(4)} - a_2h_1^{(4)})h_3\\ 
&+\left(8(h_1h_2'-h_1'h_2) + 10(a_1h_2''' - a_2h_1''')\right)h_3'\\
&+18(a_1h_2'' - a_2h_1'')h_3'' + 12(a_1h_2' - a_2h_1')h_3'''\\
&+2h_3'''\int_{S^1}h_1'h_2''dx - h_1'''\int_{S^1}h_2'h_3''dx
+h_2'''\int_{S^1}h_1'h_3''dx\\
&+2a_3(h_1'''h_2' - h_1'h_2''')
+2a_3(h_1h_2^{(4)} - h_1^{(4)}h_2)
+a_3(a_1h_2^{(6)} - a_2h_1^{(6)})
\endaligned} \\
\vphantom{\int_A^A} \\
{\aligned
&\int_{S^1} h_3'''(a_1h_2'''-a_2h_1''')dx\\ 
&+\int_{S^1}2h_3'(h_1h_2'''-h_1'''h_2 - 2h_1'h_2'' + 2h_1''h_2')dx
\endaligned}
\endpmatrix 
$$
which coincides with formula (2.3) in Misiolek \cit!{22}.
This in turn leads to the following expression for the sectional
curvature  
$\left\langle4\Cal R\left(\tbinom{h_1}{a_1},\tbinom{h_2}{a_2}\right) 
     \tbinom{h_1}{a_1},\tbinom{h_2}{a_2}\right\rangle = $ 
$$\align 
=\int_{S^1}\Bigl( 
&4(h_1h_2''-h_1''h_2)h_1h_2 +8(h_1h_2'-h_1'h_2)h_1'h_2\\ 
&+2(a_1h_2^{(4)}-a_2h_1^{(4)})h_1h_2+10(a_1h_2'''-a_2h_1''')h_1'h_2\\ 
&+18(a_1h_2''-a_2h_1'')h_1''h_2\\ 
&+12(a_1h_2'-a_2h_1')h_1'''h_2+2\om(h_1,h_2)h_1'''h_2\\ 
&-h_1'''\om(h_2,h_1)h_2+h_2'''\om(h_1,h_1)h_2\\ 
&+2(h_1'''h_2' - h_1'h_2''')a_1h_2\\ 
&+2(h_1h_2^{(4)}-h_1^{(4)}h_2)a_1h_2\\  
&+(a_1h_2^{(6)}-a_2h_1^{(6)})a_1h_2 \\ 
&+(4h_1'h_1h_2'''+2h_1h_1'h_2'''+a_1h_1'''h_2'''\\ 
&\qquad-4h_2'h_1h_1'''-2h_2h_1'h_1'''-a_2h_1'''h_1''')a_2    \Bigr)\;dx\\ 
\allowdisplaybreak 
=\int_{S^1}\Bigl( 
&-4[h_1,h_2]^2 
+4(a_1h_2-a_2h_1)(h_1h_2^{(4)}-h_1'h_2'''+h_1'''h_2'-h_1^{(4)}h_2)\\ 
&-(h_2''')^2a_1^2 +2h_1'''h_2'''a_1a_2 -(h_1''')^2a_2^2 \Bigr)\;dx\\ 
+3\om(&h_1,h_2)^2 .\\ 
\endalign$$
This formula shows that the sign of the sectional curvature is not
constant. Indeed, choosing $h_1(x) = \sin x, h_2(x) = \cos x$ we get 
$-\pi(8 + a_1^2 + a_2^2 - 3\pi)$ which can be positive and negative
by choosing the constants $a_1, a_2$ judiciously.

\subhead\nmb.{5.3}. Jacobi fields \endsubhead 
A Jacobi field $y=\binom{y}{b}$ along a geodesic with velocity field  
$\binom{u}{a}$ is a solution of the partial differential  
equation \nmb!{3.4}\thetag1 which in our case looks as follows.  
$$\align 
\binom{y_{tt}}{b_{tt}} 
&= \left[\ad\binom{y}{b}^\top+\ad\binom{y}{b}, 
                     \ad\binom{u}{a}^\top\right]\binom{u}{a} \\ 
&\quad- \ad\binom{u}{a}^\top \binom{y_t}{b_t}  
     -\al\binom{u}{a}\binom{y_t}{b_t} + \ad\binom{u}{a}\binom{y_t}{b_t}\\ 
&= \left[\pmatrix 3y_x & y_{xxx} \\  
                  \om(y) & 0 \endpmatrix, 
           \pmatrix 2u_x+u\partial_x & u_{xxx}\\ 
                     0 & 0 \endpmatrix      \right]\binom{u}{a} \\ 
&\quad+ \pmatrix -2u_x-4u\partial_x-a\partial_x^3 &
-u_{xxx}\\ 
                   \om(u) & 0 \endpmatrix         \binom{y_t}{b_t}, \\ 
\endalign$$ 
which leads to  
$$\align 
y_{tt}&=-u(4y_{tx}+3uy_{xx}+ay_{xxxx}) -u_x(2y_t+2ay_{xxx})\tag1\\ 
&\qquad-u_{xxx}(b_t+\om(y,u)-3ay_x)-ay_{txxx}, \\ 
b_{tt}&=\om(u,y_t)+\om(y,3u_xu)+\om(y,au_{xxx}).\tag2 
\endalign$$ 
Equation \thetag2 is equivalent to: 
$$ 
b_{tt} = \int_{S^1}(-y_{txxx}u + y_{xxx}(3u_xu+au_{xxx}))dx. \tag{2$'$}
$$ 
Next, let us show that the integral term in equation
\thetag1 is constant:
$$
b_t+\om(y,u) = b_t +\int_{S^1} y_{xxx}u\,dx =: B_1. \tag3 
$$
Indeed its $t$-derivative along the geodesic for $u$ (that is, $u$
satisfies the Korteweg-de Vries equation)
coincides with \thetag{2$'$}:
$$
b_{tt} + \int_{S^1} (y_{txxx}u + y_{xxx}u_t)\, dx  
= b_{tt} + \int_{S^1} (y_{txxx}u + y_{xxx}(-3u_xu-au_{xxx}))\, dx = 0. 
$$ 
Thus $b(t)$ can be explicitly solved from \thetag3 as 
$$ 
b(t) = B_0 + B_1t - \int_a^t\!\!\int_{S^1}y_{xxx}u\,dx\,dt.\tag4 
$$ 
The first component of the Jacobi equation on the Virasoro-Bott group  
is a genuine partial differential equation. Thus the Jacobi equations
are given by the following  system: 
$$\align 
y_{tt} &= -u(4y_{tx}+3uy_{xx}+ay_{xxxx}) -u_x(2y_t+2ay_{xxx})\\ 
&\qquad-u_{xxx}(B_1-3ay_x)-ay_{txxx}, \tag5\\ 
u_t&= -3u_xu - au_{xxx},\\ 
a&= \text{ constant, } 
\endalign$$ 
where $u(t,x),y(t,x)$ are either smooth functions in  
$(t,x)\in I\x S^1$ or in  
$(t,x)\in I\x \Bbb R$, where $I$ is an interval or $\Bbb R$, and  
where in the latter case $u$, $y$, $y_t$  
have compact support with respect to $x$.

Choosing $u = c \in \Bbb R$, a constant, these equations coincide with
(3.1) in Misiolek \cit!{22} where it is shown by direct inspection that
there are solutions of this equation which vanish at non-zero values
of $t$, thereby concluding that there are conjugate points along 
geodesics emanating from the identity element of the Virasoro-Bott
group on $S^1$.
 
\subhead\nmb.{5.4}. The weak symplectic structure on the space of  
Jacobi fields on the Virasoro Lie algebra \endsubhead 
Since the Korteweg - de Vries equation has local solutions  
depending smoothly on the initial conditions (and global solutions if  
$a\ne 0$), we expect that 
the space of all Jacobi fields exists and is isomorphic to the space 
of all initial data  
$(\Bbb R\x_\om \X(S^1))\x(\Bbb R\x_\om \X(S^1))$. The weak symplectic  
structure is given in section \nmb!{3.6}: 
$$\align 
\si\left(\binom{y}{b},\binom{z}{c}\right)  
&= \left\langle \binom{y}{b},\binom{z_t}{c_t}\right\rangle  
     - \left\langle \binom{y_t}{b_t},\binom{z}{c}\right\rangle 
     + \left\langle \left[\binom{u}{a},\binom{y}{b}\right], 
          \binom{z}{c}\right\rangle \\ 
&\qquad- \left\langle \binom{y}{b}, 
          \left[\binom{u}{a},\binom{z}{b}\right] \right\rangle  
     - \left\langle \left[\binom{y}{b},\binom{z}{c}\right], 
          \binom{u}{a}\right\rangle \\ 
&= \int_{S^1\text{ \!\!or }\Bbb R} (yz_t-y_tz+2u(yz_x-y_xz))\,dx \\ 
&\qquad+b(c_t+\om(z,u)) -c(b_t+\om(y,u)) -a\om(y,z) \\ 
&= \int_{S^1\text{ \!\!or }\Bbb R} (yz_t-y_tz+2u(yz_x-y_xz))\,dx \tag1\\ 
&\qquad+bC_1 -cB_1 -a\int_{S^1\text{ \!\! or }\Bbb R}y'z''\,dx,
\endalign$$ 
where the constant $C_1$ relates to $c$ as $B_1$ does to $b$, see 
\nmb!{5.3}.\thetag3 and \thetag4.


\Refs 
 
 
\ref 
\no \cit0{1} 
\by Arnold, V.I. 
\paper Sur la g\'eometrie diff\'erentielle des groupes de Lie de  
dimension infinie et ses applications \`a l'hydrodynamique des  
fluides parfaits 
\jour Ann. Inst. Fourier 
\vol 16 
\yr 1966 
\pages 319--361 
\endref 
 
\ref 
\no \cit0{2} 
\by Arnold, V.I. 
\paper An a priori estimate in the theory of hydrodynamic stability 
\paperinfo Russian 
\jour Izvestia Vyssh. Uchebn. Zaved. Matematicka 
\vol 54,5 
\yr 1966 
\pages 3-5 
\endref 
 
\ref 
\no \cit0{3} 
\by Arnold, V.I. 
\book Mathematical methods of classical mechanics 
\bookinfo Graduate Texts in Math. 60 
\publ Springer-Verlag 
\publaddr New York, Heidelberg 
\yr 1978 
\endref 
 
\ref 
\no \cit0{4} 
\by Binz, E. 
\paper Two natural metrics and their covariant derivatives on a  
manifold of embeddings 
\jour Monatsh. Math. 
\vol 89 
\yr 1980 
\pages 275--288 
\endref 
 
\ref  
\no \cit0{5} 
\by Binz, Ernst; Fischer, Hans R. 
\paper The manifold of embeddings of a closed manifold 
\inbook Proc. Differential geometric methods in 
theoretical physics, Clausthal 1978 
\publ Springer Lecture Notes in Physics 139  
\yr 1981  
\endref  
 
\ref  
\no \cit0{6} 
\by Fr\"olicher, A.; Kriegl, A. \book Linear 
spaces and differentiation theory \bookinfo Pure and Applied 
Mathematics \publ J. Wiley \publaddr Chichester \yr 1988 \endref 
 
\ref  
\no \cit0{7} 
\by Fuks, D. B.  
\book Cohomology of infinite dimensional Lie algebras  
\bookinfo (Russian)  
\publ Nauka  
\publaddr Moscow  
\yr 1984  
\transl English 
\bookinfo Contemporary Soviet Mathematics 
\publ Consultants Bureau (Plenum Press) 
\publaddr New York 
\yr 1986 
\endref 

\ref
\no \cit0{8} 
\by Gelfand, I.M.; Dorfman, I.Y.
\paper Hamiltonian operators and the algebraic structures connected 
with them
\jour Funct. Anal. Appl. 
\vol 13
\yr 1979
\pages 13--30
\endref

 
\ref 
\no \cit0{9} 
\by Kainz, Gerd 
\paper A note on the manifold of immersions 
and its Riemannian curvature 
\jour Monatshefte f\"ur Mathematik 
\vol 98  
\yr 1984 
\pages 211-217 
\endref 

\ref
\no \cit0{10}
\by Kirillov, A.A.
\paper The orbits of the group of diffeomorphisms of the circle, and 
local Lie superalgebras 
\jour Funct. Anal. Appl.
\vol 15
\yr 1981
\pages 135--136
\endref
 
\ref
\no \cit0{11}
\by Kobayashi, S.; Nomizu, K.
\book Foundations of Differential Geometry. Vol. I.\,
\publ J. Wiley-Intersci-ence 
\yr 1963 \moreref
\book Vol. II
\yr 1969
\endref

\ref   
\no \cit0{12} 
\by Kriegl, A.; Michor, P. W.   
\paper A convenient setting for real analytic mappings  
\jour Acta Mathematica   
\vol 165  
\pages 105--159  
\yr 1990   
\endref 
 
\ref 
\no \cit0{13} 
\by Kriegl, A.; Michor, P. W. 
\paper Regular infinite dimensional Lie groups 
\jour J. Lie Theory (http://ww w.emis.de/journals/JLT) 
\vol 7,1 
\yr 1997  
\pages  
\finalinfo ESI Preprint 200 
\endref 

\ref  
\no \cit0{14} 
\by Kriegl, A.; Michor, P. W.  
\book The Convenient Setting for Global Analysis   
\bookinfo Surveys and Monographs 53
\publ AMS 
\publaddr Providence 
\yr 1997 
\endref 
 
\ref 
\no \cit0{15} 
\by Marsden, J; Ratiu, T. 
\book Introduction to mechanics and symmetry 
\publ Springer-Verlag 
\publaddr New York, Berlin, Heidelberg 
\yr 1994 
\endref 
 
\ref  
\no \cit0{16} 
\by Michor, P. W. 
\paper  Manifolds of smooth maps  
\jour Cahiers Topol. Geo. Diff. 
\vol 19  
\yr 1978 
\pages 47--78 
\endref 
 
\ref  
\no \cit0{17} 
\by Michor, P. W. 
\paper Manifolds of smooth maps II: The Lie group 
of diffeomorphisms of a non compact smooth manifold 
\jour Cahiers Topol. Geo. Diff.  
\vol 21  
\yr 1980 
\pages 63--86 
\endref 
 
\ref  
\no \cit0{18} 
\by Michor, P. W. 
\paper Manifolds of smooth maps III: The principal bundle 
of embeddings of a non compact smooth manifold 
\jour Cahiers Topol. Geo. Diff.  
\vol 21  
\yr 1980 
\pages 325--337 
\endref 
 
 
\ref 
\no \cit0{19} 
\by Michor, P. W. \book Manifolds of differentiable mappings 
\publ Shiva \yr 1980 \publaddr Orpington \endref 
 
\ref  
\no \cit0{20} 
\by Michor, P. W. \paper Manifolds of smooth 
mappings IV: Theorem of De~Rham \jour Cahiers Top. Geo. Diff.  
\vol 24 \yr 1983 \pages 57--86 \endref 
 
\ref  
\no \cit0{21} 
\by Michor, P. W. \paper Gauge theory for 
diffeomorphism groups \inbook Proceedings of the Conference on 
Differential Geometric Methods in Theoretical Physics, Como 
1987, K. Bleuler and M. Werner (eds.) 
\publ Kluwer \publaddr Dordrecht \yr 1988 \pages 345--371 \endref

\ref  
\no \cit0{22} 
\by Misiolek, G. \paper Conjugate points in the Bott-Virasoro 
group and the KdV equation \jour
Proc. Amer. Math. Soc.  
\vol 125 \yr 1997 \pages 935--940 \endref 

\ref
\no \cit0{23} 
\by Ovsienko, V.Y; Khesin, B.A.
\paper Korteweg--de~Vries superequations as an Euler equation
\jour Funct. Anal. Appl.
\vol 21
\yr 1987
\pages 329--331
\endref

\ref
\no \cit0{24} 
\by Segal, G.
\paper The geometry of the KdV equation
\jour Int. J. Mod. Phys. A
\vol 6
\yr 1991
\pages 2859--2869
\endref
  
\ref  
\no \cit0{25} 
\by Weinstein, Alan 
\paper Symplectic manifolds and their Lagrangian manifolds 
\jour Advances in Math. 
\vol 6 
\pages 329--345 
\yr 1971 
\endref 
 
\endRefs
\enddocument